\documentclass[11pt]{article}

\usepackage{latexsym,color,graphicx,amsmath,amssymb}

\def\reals{{\mathbb R}}
\def\cplx{{\mathbb C}}
\def\eps{{\varepsilon}}

\newtheorem{theorem}{Theorem}[section]

\begin{document}

\title{Incidences between points on a variety and planes in $\reals^3$\thanks{%
Work on this paper by Noam Solomon and Micha Sharir was
supported by Grant 892/13 from the Israel Science Foundation.
Work by Micha Sharir was also supported
by Grant 2012/229 from the U.S.--Israel Binational Science Foundation,
by the Blavatnik Research Fund for Computer Science at Tel Aviv University,
by the Israeli Centers of Research Excellence (I-CORE) program (Center No.~4/11), and
by the Hermann Minkowski-MINERVA Center for Geometry
at Tel Aviv University.
}}

\author{
Micha Sharir\thanks{%
School of Computer Science, Tel Aviv University,
Tel Aviv 69978, Israel.
{\sl michas@post.tau.ac.il} }
\and
Noam Solomon\thanks{%
School of Computer Science, Tel Aviv University,
Tel Aviv 69978, Israel.
{\sl noam.solom@gmail.com} }
}

\maketitle

\begin{abstract}
In this paper we establish an improved bound for the number of incidences between
a set $P$ of $m$ points and a set $H$ of $n$ planes in $\reals^3$, provided that
the points lie on a two-dimensional nonlinear irreducible algebraic variety $V$ of constant degree.
Specifically, the bound is
$$
O\left( m^{2/3}n^{2/3} + m^{6/11}n^{9/11}\log^\beta(m^3/n) + m + n +
\sum_\ell |P_\ell|\cdot |H_\ell| \right) ,
$$
where the constant of proportionality and the constant exponent $\beta$ depend on the
degree of $V$, and where the sum ranges over all lines $\ell$ that are fully contained in
$V$ and contain at least one point of $P$, so that, for each such $\ell$,
$P_\ell = P\cap\ell$ and $H_\ell$ is the set of the planes of $H$ that contain $\ell$.
In addition, $\sum_\ell |P_\ell| = O(m)$ and $\sum_\ell |H_\ell| = O(n)$.

This improves, for this special case, the earlier more general bound of Apfelbaum and Sharir~\cite{ApS}
(see also \cite{BK,ET05}).

This is a generalization of the incidence bound for points and circles in the
plane~\cite{ANP+,ArS,MT}, and is based on a recent result of Sharir and Zahl~\cite{SZ}
on the number of cuts that turn a collection of algebraic curves into pseudo-segments.

The case where $V$ is a quadric is simpler to analyze, does not require the result of \cite{SZ},
and yields the same bound as above, with $\beta=2/11$.

We present an interesting application of our results to a problem, studied by
Rudnev~\cite{Rud}, on obtaining a lower bound on the number of distinct cross-ratios
determined by $n$ real points, where our bound leads to a slight improvement in Rudnev's bound.
\end{abstract}

\section{Introduction}

Let $P$ be a set of $m$ points, and $H$ a set of $n$ planes in
$\reals^3$. Assume that $P$ is contained in some two-dimensional nonlinear and irreducible
algebraic variety $V$ of constant degree $D$. We wish to bound the size $I(P,H)$ of
the \emph{incidence graph} $G(P,H)$, whose edges connect all pairs $(p,h)\in P\times H$
such that $p$ is incident to $h$. In general, and in this special setup too, $I(P,H)$ can be
as large as the maximum possible value $|P|\cdot|H|$, by placing all the points of $P$
on a line, and make all the planes of $H$ contain this line, in which case $G(P,H) = P\times H$.
The bound that we are going to obtain will of course acknowledge this possibility,
and will be obtained by decomposing $G(P,H)$ as the disjont union
$G_0(P,H) \cup \bigcup_i P_i\times H_i$, where, for each $i$, there exists a line
$\ell_i$ fully contained in $V$, so that $P_i = P\cap\ell_i$ and $H_i$ is the set
of planes of $H$ that contain $\ell_i$. The residue graph $G_0(P,H)$ represents the
``accidental'' incidences, those that cannot be ``explained'' in terms of lines that
are contained in $V$ and induce large complete bipartite subgraphs of $G(P,H)$, and we
provide an explicit bound for the number $I_0(P,H) = |G_0(P,H)|$ of these incidences.
The quality of the bound will be measured by $I_0(P,H)$ and by $\sum_i |P_i|$ and $\sum_i |H_i|$.

Earlier works on point-plane incidences have considered the general
setup, where the points of $P$ are arbitrarily placed in $\reals^3$.
Initial partial results have been obtained by Edelsbrunner, Guibas
and Sharir~\cite{EGS}. More recently, Apfelbaum and
Sharir~\cite{ApS} (see also Brass and Knauer~\cite{BK} and Elekes
and T{\'{o}}th~\cite{ET05}) have shown that if the incidence graph
does not contain a copy of $K_{r,s}$, for constant parameters $r$
and $s$, then $I(P,H)=O(m^{3/4}n^{3/4}+m+n)$. In more generality,
Apfelbaum and Sharir~\cite{ApS} have shown that if $I=I(P,H)$ is
significantly larger than this bound, then $G(P,H)$ contains a large
complete bipartite subgraph $P'\times H'$, such that $|P'|\cdot |H'|
= \Omega(I^2/(mn)) - O(m+n)$. Moreover, as also shown in~\cite{ApS}
(slightly improving a similar result of Brass and Knauer~\cite{BK}),
$G(P,H)$ can be expressed as the union of complete bipartite graphs
$P_i\times H_i$ so that $\sum_i (|P_i|+|H_i|) =
O(m^{3/4}n^{3/4}+m+n)$. (This is a specialization to the case $d=3$
of a similar result of \cite{ApS,BK} in any dimension $d$.) A more
recent result of Zahl~\cite{Za} gives a grand generalization of
these bounds, to incidences between points and constant-degree
algebraic surfaces for which the incidence graph does not contain
$K_{r,3}$ as a subgraph; that is, in situations where every triple
of the given surfaces intersect in a constant number of points.
Zahl's general bound is $O(m^{3/4}n^{3/4}+m+n)$, as above. We note
that Fox et al.~\cite{FPSSZ} present a more general framework that
includes incidences problems of many kinds, and yields, for the case
under consideration, almost the same bound, namely $O(m^{3/4
+\eps}n^{3/4 + \eps} + m + n)$, for any $\eps$, where the constant
of proportionality depends on $\eps$.

Nevertheless, as we show in this paper, the bound can be substantially improved
when all the points of $P$ lie on a constant-degree variety $V$. Our main result
is the following theorem.
\begin{theorem} \label{main}
For a set $P$ of $m$ points on some two-dimensional, nonlinear,
non-conical\footnote{That is, $V$ is not the union of lines through
a common point.} and irreducible algebraic variety $V$ of constant
degree $D$ in $\reals^3$, and for a set $H$ of $n$ planes in
$\reals^3$, the incidence graph $G(P,H)$ can be decomposed as
\begin{equation} \label{gph}
G(P,H) = G_0(P,H) \cup \bigcup_i (P_i\times H_i) ,
\end{equation}
such that, for each $i$, there exists a line $\ell_i \subset V$ so that
$P_i = P\cap \ell_i$ and $H_i$ is the subset of planes of $H$ that contain $\ell_i$, and we have
\begin{align} \label{eq:main}
|G_0(P,H)| & = O\left( m^{2/3}n^{2/3} + m^{6/11}n^{9/11}\log^\beta(m^3/n) + m + n \right) , \quad\text{and} \\
\sum_i |P_i| & = O(m) \quad\text{and}\quad \sum_i |H_i| = O(n) \nonumber ,
\end{align}
where the constant exponent $\beta$ and the constants of proportionality depend on the degree $D$ of $V$.

If $V$ is a (non-linear) cone, the same result holds if we exclude
the apex of $V$ from $P$; the apex adds at most $O(n)$ incidences.
\end{theorem}

We remark that the problem studied in this paper is not as specialized as it
might appear. It is a generalization of the problem of
bounding the number of incidences between points and circles of arbitrary radii
in the plane. By applying the standard lifting transformation, which maps each
point $(x,y)$ in the plane to the point $(x,y,x^2+y^2)$ on the paraboloid
$z=x^2+y^2$ in 3-space, and maps each circle $(x-a)^2+(y-b)^2=r^2$ to the plane
$z=2ax+2by+(r^2-a^2-b^2)$, the problem becomes that of bounding the number
of incidences between points on the paraboloid and arbitrary planes in three dimensions.
This transformation is used in the second step of the analysis of the point-circle
incidence problem, as in \cite{ANP+,ArS}, and the bound derived there is very similar
to the one in (\ref{eq:main}), except that, in the case of circles, $G(P,H)$ does not
contain large complete bipartite graphs (because the paraboloid does not contain any line),
and the exponent $\beta$ is equal to the sharper value $2/11$; see~\cite{ANP+,ArS,MT} for details.

As a matter of fact, we extend this observation in Section~\ref{sec:quadric},
and derive such an improved variant of Theorem~\ref{main} for any quadric $V$.
The resulting bound is the same, but the analysis is simplified, and the exponent
$\beta$ in (\ref{eq:main}) is still $2/11$; see Theorem~\ref{thm:quadric} for the precise formulation.

As an application, we note that our bound leads to a slight improvement in a recent lower bound,
due to Rudnev~\cite{Rud}, on the number of distinct cross-ratios determined by $n$ real numbers,
improving his bound from $\Omega(n^{24/11}/\log^{6/11}n)$ to $\Omega(n^{24/11})$. More details
are provided in Section~\ref{sec:cross}.

In closing the introduction, we note that the terms $P_i\times H_i$ can only arise when
$V$ fully contain the corresponding line $\ell_i$. This issue, of analyzing lines fully contained
in algebraic varieties in $\reals^3$ (or in $\cplx^3$) has received considerable attention
in recent years, starting with Guth and Katz's work~\cite{GK2}. Briefly, if $V$ is a
\emph{ruled surface},\footnote{%
  That is, each point on $V$ is incident to a line that is fully contained in $V$.}
it contains infinitely many lines, and then the set of lines that form the complete
bipartite decompositions might be large (up to $\Theta(\min\{m,n\})$ lines). That is,
ruled surfaces might yield a large number of complete bipartite subgraphs in $G(P,H)$
(but then these graphs are smaller on average). In contrast, as follows from the
Monge--Cayley--Salmon theorem (see, e.g., Guth and Katz~\cite{GK2}), a non-ruled
irreducible variety of degree $D$ can fully contain at most $D(11D-24)$ lines. Hence,
if $V$ is not ruled, the decomposition in (\ref{gph}) can consist of only a constant
number of complete bipartite graphs.

\section{Proof of Theorem~\ref{main}}

Let $P$, $V$, $H$, $m$, and $n$ be as above. To obtain the representation (\ref{gph})
and the corresponding bounds (\ref{eq:main}), we first derive a weaker bound, and then
improve it via a suitable decomposition of dual space, similar to the way it has been
done for circles in \cite{ANP+,ArS}, and in more generality in \cite{SZ}.

\paragraph{A weak basic bound.}
Let $L$ denote the set of lines fully contained in $V$ and containing
at least one point of $P$. For each line $\ell\in L$ we form the
bipartite subgraph $P_\ell \times H_\ell$ of $G(P,H)$, where $P_\ell = P\cap\ell$
and $H_\ell$ is the set of all the planes of $H$ that contain $\ell$.
(Recall the discussion at the end of the introduction: There are at most $O(1)$
lines $\ell$ (and corresponding complete bipartite graphs) when $V$ is not ruled.)
As shown, e.g., in Sharir and Solomon~\cite[Lemma 5]{SS3d}, except for at
most one \emph{exceptional} point, each point $p \in V$ is incident
to at most $D^2=O(1)$ lines of $L$.

If $V$ contains an exceptional point $p_0$ then it must be a cone with $p_0$ as
an apex, and (since $V$ is nonlinear) all the lines contained in $V$ pass through $p_0$.
In this case we remove $p_0$ from $P$, losing at most $n$ incidences,
and re-apply the analysis to the set $P\setminus\{p_0\}$, which we continue to denote by $P$.

Ignoring then the exceptional point, if any, the preceding property
implies that $\sum_{\ell} |P_\ell|=O(m)$. (That is, even when $V$ is ruled,
only finitely many lines $\ell$ can arise.) As $V$ is irreducible and
nonlinear, it does not contain any of the planes in $H$. Thus, for
each $h\in H$, the intersection $h\cap V$ is a (plane) algebraic curve of
degree at most $D$, and can therefore contain at most $D=O(1)$ lines
of $L$ (as follows, e.g., from the generalized version of B\'ezout's
theorem~\cite{Fu84}, or by other, more direct arguments).
This implies that $\sum_{\ell} |H_\ell|=O(n)$. Retaining only the lines that
contain at least one point of $P$ and are contained in at least one plane of $H$,
we obtain a finite decomposition $\bigcup_\ell \left( P_\ell\times H_\ell \right)$
as a portion of $G(P,H)$, and $\sum_\ell |P_\ell| = O(m)$, and $ \sum_\ell |H_\ell| = O(n)$.

For each $h\in H$, put $\gamma_h := (h\cap V) \setminus \bigcup L$.
As just noted, each $\gamma_h$ is at most one-dimensional. By
construction, it does not contain any line, and it might also be
empty (for this or for other reasons). Note that if $h\cap V$ does
contain a line $\ell$, then $\ell\in L$, and the incidences between
$h$ and the points of $P$ on $\ell$ are all already recorded in
$P_\ell\times H_\ell$. Finally, we ignore the isolated points of
$\gamma_h$. By Harnack's curve theorem~\cite{Har}, the number of
such points on a plane curve of degree at most $D$ is $O(D^2) = O(1)$,
so the isolated points contribute a total of at most $O(n)$ incidences.
Let $G_0(P,H)$ denote the remaining portion of $G(P,H)$, after
pruning away the complete bipartite graphs $P_\ell\times H_\ell$.
This gives us the decomposition in (\ref{gph}). We further remove
from $G_0(P,H)$ all the $O(n)$ incidences involving isolated points
on the curves $\gamma_h$, and the $O(n)$ incidences with the
exceptional point of $V$, if it exists and belongs to $P$, continue
to denote the resulting graph as $G_0(P,H)$, for notational
convenience, and put $I_0(P,H) = |G_0(P,H)|$.

Let $\Gamma$ denote the set of the $n$ curves $\gamma_h$, for $h\in H$. The curves of $\Gamma$
are algebraic curves of degree at most $D$, and any pair of curves
$\gamma_h$, $\gamma_{h'}\in \Gamma$ intersect in at most $D=O(1)$ points.
Indeed, any of these points is an intersection point of $V$ with the line
$\ell=h\cap h'$; if $\ell$ is contained in $V$ then, by construction, it
is removed from both curves, and if $\ell$ is not contained in $V$, it can
meet it in at most $D$ points. In particular, the curves $\gamma_h$
are distinct and non-overlapping.

Note that $I_0(P,H)$ is equal to the number of incidences $I(P,\Gamma)$
between the points of $P$ and the curves of $\Gamma$. To bound the latter
quantity we proceed as follows.

We slightly tilt the coordinate frame to make it generic, and then project the
curves of $\Gamma$ onto the $xy$-plane. A suitable choice of the tilting guarantees
that (i) no pair of intersection points or points of $P$ project to the same point;
(ii) if $p$ is not incident to $\gamma_h$ then the projections of $p$ and of
$\gamma_h$ remain non-incident; and (iii) no pair of curves in $\Gamma$ have
overlapping projections. In addition, by
construction, no curve of $\Gamma$ contains any (vertical) segment. Let $P^*$ and
$\Gamma^*$ denote, respectively, the set of projected points and the set of
projected curves; the latter is a set of $n$ plane algebraic curves of constant maximum
degree $D$ (since the projection of a plane curve onto another plane clearly cannot increase its degree).
Moreover, $I(P,\Gamma)$ is equal to the number of incidences between $P^*$ and $\Gamma^*$.

By the result of Sharir and Zahl~\cite{SZ}, applied to $\Gamma^*$,
the curves of $\Gamma^*$ can be cut into $O(n^{3/2}\log^\kappa n)$ connected Jordan subarcs,
where the constant exponent $\kappa$ and the constant of proportionality depend on $D$,
so that each pair of subarcs intersect at most once. Using standard terminology,
the subarcs form a collection of \emph{pseudo-segments}.

We can now apply the crossing-lemma technique of Sz\'ekely~\cite{Sze}, exactly as
was done in Sharir and Zahl~\cite{SZ} or in earlier works. Since the resulting subarcs
form a collection of pseudo-segments, and the number of their intersections is $O(n^2)$,
Sz\'ekely's analysis yields the bound
\begin{equation} \label{i0ph}
I_0(P,H) = I(P,\Gamma) = I(P^*,\Gamma^*) =
O\left( m^{2/3}n^{2/3} + m + n^{3/2}\log^\kappa n \right).
\end{equation}

Adding incidences recorded in the complete bipartite decomposition, as constructed
above, and those involving isolated points and the exceptional point, if any,
we get our initial (weak) bound.
\begin{equation} \label{weakinc}
I(P,H) = O\left( m^{2/3}n^{2/3} + m + n^{3/2}\log^\kappa n + \sum_{\ell} |P_{\ell}| \cdot |H_{\ell}| \right) ,
\end{equation}
where $\bigcup_{\ell} \left( P_{\ell}\times H_{\ell} \right)$ is contained in the incidence
graph $G(P,H)$, is a finite union over lines $\ell\subset V$,
and $\sum_{\ell} |P_{\ell}| = O(m)$ and $\sum_{\ell} |H_{\ell}| = O(n)$.

\paragraph{The case $m=O(n^{1/3})$.}
Before proceeding to improve the bound in (\ref{weakinc}), we first
dispose of the case $m=O(n^{1/3})$. As above, we first remove all
the complete bipartite graphs $P_\ell\times H_\ell$, for $\ell\in L$,
from $G(P,H)$.
We then proceed to estimate $I_0(P,H)$ as follows. We call a
plane $h\in H$ \emph{strongly degenerate} (or \emph{degenerate} for short)
if all the points of $P\cap h$ are collinear. We claim that the number
of incidences between points of $P$ and non-degenerate planes is
$O(m^3+n)=O(n)$. Indeed, first discard the planes $h\in H$
containing at most two points of $P$, losing at most $O(n)$ incidences.
For an incidence between a point $p\in P$ and a surviving non-degenerate
plane $h\in H$, there exist (at least) two distinct points
$q,s \in (P\setminus \{p\})\cap h$ that are not collinear with $p$.
The ordered triple $(p,q,s)$ therefore \emph{uniquely} accounts for the
incidence between $p$ and $h$, and there are $O(m^3)=O(n)$ such triples.

To bound the number of incidences between the
points of $P$ and the degenerate planes of $H$, fix a degenerate
plane $h\in H$, and assume that $m_h:=|h\cap P| \ge D+1$; the
overall number of incidences on the other planes is at most $Dn = O(n)$.
By assumption, all points of $h\cap P$ lie on a common line
$\ell$. Since $\ell$ contains at least $D+1$ points, it must be
contained in $V$, so the incidences involving $h$ are all recorded
in the complete bipartite graph $P_\ell \times H_\ell$.
In other words, we have shown, for $m = O(n^{1/3})$,
\begin{equation} \label{smallm}
I(P,H) = O\left( n + \sum_\ell |P_\ell| \cdot |H_\ell| \right) ,
\end{equation}
where, as above, $\bigcup_\ell \left( P_\ell\times H_\ell \right)$ is contained
in the incidence graph $G(P,H)$,
and $\sum_{\ell} |P_{\ell}| = O(m)$ and $\sum_{\ell} |H_{\ell}| = O(n)$.
In the rest of the analysis, we thus assume that $m=\Omega(n^{1/3})$.

\paragraph{Improving the bound.}
As in the analysis of incidences between points and circles (or pseudo-circles) in
\cite{ANP+,ArS}, and the more general treatment in~\cite{SZ}, the first two terms
in (\ref{i0ph}) dominate when $m = \Omega(n^{5/4}\log^{3\kappa/2} n)$. When $m$
is smaller, the third term, which is independent of $m$, is the one that dominates,
and we then sharpen it as follows (the same general approach is also used in
\cite{ANP+,ArS,SZ}).

We apply any standard duality to 3-space, which maps points to planes and planes to points,
and preserves point-plane incidences. We denote by $a^*$ the dual image (plane or point)
of an object $a$ (point or plane). We thus get a set $H^*$ of $n$ points, and a set $P^*$
of $m$ planes in dual 3-space.\footnote{%
  The planes are special, because they come from points on a variety, but we do not know how
  to exploit this property.}
We choose a parameter $r$, to be fixed later, and construct a
\emph{$(1/r)$-cutting} $\Xi$ for $P^*$ (see, e.g., \cite{Chaz}), which partitions the dual
space into $O(r^3)$ simplices, each crossed by at most $m/r$ dual planes.

\paragraph{Incidences with boundary dual points.}
Let us first bound the number of incidences between the dual points
that lie on the boundaries of the simplices and the dual planes.

\smallskip

\noindent
(i) Each dual point that lies in the relative interior of a facet of some
simplex of $\Xi$ (but not on any edge of another simplex) has one incidence with
the dual plane that contains that facet (if any), for a total of $O(n)$ such incidences.
If such a point $h^*$ is incident to another dual plane $p^*$, then
$p^*$ must cross both simplices adjacent to the facet containing $h^*$.
We then assign $h^*$ to one of these two simplices, and the relevant incidence
will then be counted in the subproblem associated with that simplex
(see below).

\smallskip

\noindent
(ii) Any incidence between a dual plane $p^*$ and a dual point $h^*$
that is a vertex of some simplices of $\Xi$ can be charged to the
crossing of $p^*$ with some simplex $\sigma$ of $\Xi$ that is adjacent
to $h^*$. It follows that the number of incidences of this kind is at
most $O(r^3)\cdot m/r = O(mr^2)$. The only case where we cannot apply
this charging is when $p^*$ supports facets of all incidence simplices,
but then the incidence is counted in (i).

\smallskip

\noindent
(iii) The remaining incidences involve points that lie on edges of the
simplices; there are $O(r^3)$ such edges. Fix one of these edges $e$.
An incidence between a point $h^*\in e$ and a plane $p^*$ that crosses
$e$ (intersects $e$ but does not contain it) can be charged to the crossing
of $p^*$ with some adjacent simplex (at most four such charges can be made to
any plane-simplex pair), and, as in (ii), the number of such incidences is only $O(mr^2)$.

\smallskip

\noindent
(iv) We are thus left
with incidences between dual points on $e$ and dual planes that
fully contain $e$. Let $\ell^*$ be the line supporting $e$. Back in
the primal space, there exists a ``pre-image'' $\ell$ of $\ell^*$
such that all the dual planes $p^*$ that contain $\ell^*$ come from
points $p$ that lie on $\ell$, and all the dual points $h^*$ that
lie on $\ell^*$ come from planes $h$ that contain $\ell$. If
$\ell\subset V$ then all these incidences are recorded in the
corresponding complete bipartite graph $P_\ell\times H_\ell$
that we have already constructed and disposed of. Otherwise, $\ell$ can contain
at most $D$ points of $P$, meaning that in dual space each point on
any of these edges $e$ is incident to at most $D$ dual planes that fully
contain $e$. Since the sets of dual points lying on the (relative interiors of the) edges $e$ are pairwise
disjoint, we get a total of at most $O(nD)=O(n)$ incidences of this latter kind.

In total, the number of incidences with boundary points is $O(mr^2+n)$.

\paragraph{Incidences within the simplices.}
We now proceed to consider incidences between dual points in the
interior of the simplices of $\Xi$. For each simplex $\sigma$ of $\Xi$, let
$n_\sigma$ denote the number of points of $H^*$ in the interior of $\sigma$,
including the points that $\sigma$ has ``inherited'' from points on its facets, in
step (i) of the preceding analysis. By duplicating cells $\sigma$
for which $n_\sigma > n/r^3$, so that in each copy we take at most $n/r^3$
of these dual points (but retain all the crossing dual planes), we obtain
a collection of $O(r^3)$ simplices, each of which is crossed by at most $m/r$
dual planes and contains at most $n/r^3$ dual points; we denote the actual
number of these planes and points as $m_\sigma$ and $n_\sigma$ (the latter
notation is slightly abused, as it now refers only to a single copy (subproblem)
of $\sigma$), respectively, for each simplex-instance $\sigma$.

For each cell $\sigma$, we apply the bound (\ref{weakinc})
to the subset $P^{(\sigma)}$ of the points of $P$ whose dual planes
cross $\sigma$, and to the subset $H^{(\sigma)}$ of the planes whose dual points
lie in $\sigma$, and note that the case $m/r = O((n/r^3)^{1/3})$ does not arise,
because then we would also have $m=O(n^{1/3})$, which we have assumed not to be
the case. That is, we get, for each $\sigma$,
$$
I(P^{(\sigma)},H^{(\sigma)}) = O\left( m_\sigma^{2/3}n_\sigma^{2/3} + m_\sigma + n_\sigma^{3/2}\log^\kappa n_\sigma
 + \sum_\ell |P^{(\sigma)}_\ell| \cdot |H^{(\sigma)}_\ell| \right) ,
$$
for a suitable complete bipartite decomposition
$\bigcup_\ell \left( P^{(\sigma)}_\ell \times H^{(\sigma)}_\ell \right)$.

We sum these bounds, over the simplices $\sigma$ of $\Xi$. We note that the same line
$\ell$ may arise in many complete bipartite graphs
$P^{(\sigma)}_\ell \times H^{(\sigma)}_\ell$, but (i) all these graphs are contained in
$P_\ell \times H_\ell$, and (ii) they are edge disjoint, because each dual point $h^*$ lies
in at most one simplex. This allows us to replace all the partial subgraphs
$P^{(\sigma)}_\ell \times H^{(\sigma)}_\ell$ by the single graph
$P_\ell \times H_\ell$, for each line $\ell$. We thus get

\begin{align*}
I(P,H) & = O\left( \sum_\sigma \Big( m_\sigma^{2/3}n_\sigma^{2/3} +
 m_\sigma + n_\sigma^{3/2}\log^\kappa n_\sigma \Big)
 + mr^2 + n + \sum_\ell |P_\ell| \cdot |H_\ell| \right) \\
& = O\left( r^3 \Big( (m/r)^{2/3}(n/r^3)^{2/3} + (n/r^3)^{3/2}\log^\kappa (n/r^3) \Big)
 + mr^2 + n + \sum_\ell |P_\ell| \cdot |H_\ell| \right) \\
& = O\left( m^{2/3}n^{2/3}r^{1/3} + \frac{n^{3/2}}{r^{3/2}}\log^\kappa (n/r^3)
 + mr^2 + n + \sum_\ell |P_\ell| \cdot |H_\ell| \right) .
\end{align*}

We now choose ${\displaystyle r = \frac{n^{5/11}\log^{6\kappa/11} (m^3/n)}{m^{4/11}}}$, to
equalize (asymptotically) the first two terms in the bound, which become
$O(m^{6/11}n^{9/11}\log^{2\kappa/11}(m^3/n))$.
(Note that $1\le r\le m$ for $n^{1/3}\le m\le n^{5/4}\log^{3\kappa/2} n$.)
The third term becomes
$mr^2 = m^{3/11}n^{10/11}\log^{12\kappa/11} (m^3/n)$, which is dominated by the preceding
bound for $m=\Omega(n^{1/3})$, as is easily checked.
The complementary case $m=O(n^{1/3})$ has been handled by (\ref{smallm}).

This completes the proof of Theorem~\ref{main}.
$\Box$

\medskip

\noindent{\bf Remarks. (i)}
In retrospect, once we have reduced the problem to that of bounding the number $I(P^*,\Gamma^*$)
on incidences between the projected points and curves on the $xy$-plane, we could have applied,
as a black-box, the analysis of Sharir and Zahl~\cite{SZ}, and get a slightly weaker bound with
an additional (arbitrarily small) $\eps$ in the exponents. Exploiting the fact that we are
dealing here with planes allowed us to obtain the finer bound in (\ref{eq:main}).

\medskip

\noindent{\bf (ii)}
Note also that the preceding analysis closely follows the one in \cite{ANP+,ArS}, or,
more precisely, its enhanced version that exploits the bound $O(n^{3/2}\log n)$
of Marcus and Tardos~\cite{MT} on the number of cuts needed to turn a collection of
circles into a collection of pseudo-segments (so $\kappa=1$ for the case of circles).

\subsection{The case of a quadric} \label{sec:quadric}

The special case where $V$ is a quadric admits a simpler analysis with a slightly
improved bound. In this case the curves $\gamma_h$, as defined earlier, are quadratic,
and by construction, do not contain any line. Moreover, each pair of curves
$\gamma_h$, $\gamma_{h'}$ intersect at most twice, at the at most two intesrection points
of the line $h\cap h'$ with $V$. Hence, these curves form a family of pseudo-circles\footnote{%
  Some care has to be exercised here, because the analysis in \cite{ANP+} handles
  pseudo-circles in the plane and not on a variety in 3-space. This extension, though,
  can be handled in a rather routine manner, e.g., by breaking $V$ into a constant
  number of $xy$-monotone pieces, projecting the configuration within each piece
  onto the $xy$-plane, and applying the analysis in \cite{ANP+} to each projection
  separately. We omit here the fairly easy details.}
with a $3$-parameter algebraic representation, in the terminology of \cite{ANP+};
they inherit this representation from the $3$-parameter algebraic representation
of the planes $h$ that define them. We can therefore apply the bound in \cite[Theorem 6.6]{ANP+},
enhanced by the improvement in \cite{MT}, and obtain the following improvement.
\begin{theorem} \label{thm:quadric}
For a set $P$ of $m$ points on some two-dimensional quadric $V$ in $\reals^3$, and for a set $H$
of $n$ planes in $\reals^3$, the incidence graph $G(P,H)$ can be decomposed as
\begin{equation} \label{gphq}
G(P,H) = G_0(P,H) \cup \bigcup_i (P_i\times H_i) ,
\end{equation}
such that, for each $i$, there exists a line $\ell_i\subset V$, so that
$P_i = P\cap \ell_i$ and $H_i$ is the subset of planes of $H$ that contain $\ell_i$, and we have
\begin{align} \label{eq:mainq}
|G_0(P,H)| & = O\left( m^{2/3}n^{2/3} + m^{6/11}n^{9/11}\log^{2/11}(m^3/n) + m + n \right) , \quad\text{and} \\
\sum_i |P_i| & = O(m) \quad\text{and}\quad \sum_i |H_i| = O(n) \nonumber .
\end{align}
\end{theorem}

\subsection{An application to distinct cross-ratios} \label{sec:cross}

In a recent work~\cite{Rud}, Rudnev has obtained a superquadratic lower bound on the number of distinct
cross-ratios determined by any set $A$ of $n$ real (or complex) numbers. One of the key steps in
his analysis (see \cite[Theorem 5(i)]{Rud}) is an upper bound on the number of $k$-rich planes in
$\reals^3$, with respect to a set $P$ of $m$ points on the hyperbolic paraboloid $\Pi:\; z=xy$;
that is, planes that contain at least $k$ points of $P$. In the interest of (partial) completeness,
here is a sketch of the reduction, adapted from \cite{Rud}, that leads to this problem. In what follows,
we only consider the real case.

Let $A$ be a set of $n$ real numbers, and consider the set of M\"obius transformations on $\reals$,
each of which can be written as
$$
y = \frac{\alpha x + \beta}{\gamma x + \delta} .
$$
To simplify things, and almost without loss of generality, assume that $\gamma = 1$
(we lose transformations with $\gamma=0$, but they are trivial to handle),
so we consider transformations of the form
$$
y = \frac{\alpha x + \beta}{x + \delta} .
$$
We also assume that $\beta \ne \alpha\delta$, so as to rule out constant transformations
from the analysis.

We are given a parameter $k$, and want to bound the number $N = N_{\ge k}$ of
\emph{$k$-rich transformations}, namely those that map at least $k$ points of
$A$ to $k$ other points in $A$.

For each pair of points $(a,b)\in A\times A$, the set of transformations that map $a$ to $b$ is
$$
M(a,b) = \left\{ (\alpha,\beta,\delta) \mid ab -\alpha a + \delta b - \beta = 0 \right\} .
$$
We can represent each triple $(\alpha,\beta,\delta)$ (i.e., each M\"obius tranformation) by the
plane $z -\alpha x + \delta y - \beta = 0$. Then $M(a,b)$ is the set of all planes that pass through
the point $(a,b,ab)$ on the hyperbolic paraboloid $\Pi$. A transformation is $k$-rich if its
associated plane passes through at least $k$ such points.

This leads to the following incidence question. We have the set
$P = \{ (a,b,ab) \mid (a,b)\in A\times A\}$ of $m=n^2$ points on $\Pi$, and some set $H$ of $N$ planes
of the above form, and wish to bound $I(P,H)$, the number of incidences between $P$ and $H$.

To apply Theorem~\ref{thm:quadric} to this scenario, we note that no plane in $H$ intersects $\Pi$
in a line. To see this, we note that, for each $h\in H$, given by an equation of the form
$z -\alpha x + \delta y - \beta = 0$, for a suitable triple $(\alpha,\beta,\delta)$ of coefficients,
its intersection with $\Pi$, when projected down to the $xy$-plane, is the hyperbola
$xy -\alpha x + \delta y - \beta = 0$. This hyperbola degenerates into a pair of lines
only when $\beta=\alpha\delta$, and since we have ruled out this situation, the claim follows.

It follows that, in the notation of Theorem~\ref{thm:quadric}, we have
$$
I(P,H) = |G_0(P,H)| = O\left( m^{2/3}N^{2/3} + m^{6/11}N^{9/11}\log^{2/11}(m^3/N) + m + n \right) .
$$
We can now obtain an upper bound on the number $N_{\ge k}$ of $k$-rich transformations,
by taking $H$ to be the set of planes that correspond to these transformations, and by comparing
this bound with the lower bound $kN_{\ge k}$. That is, we have
$$
kN_{\ge k} = O\left( m^{2/3}N_{\ge k}^{2/3} + m^{6/11}N_{\ge k}^{9/11}\log^{2/11}(m^3/N_{\ge k}) + m + N_{\ge k} \right) ,
$$
from which we easily obtain, substituting $m=n^2$ and assuming that $k$ is at least some
suitable sufficiently large constant $k_0$,
\begin{equation} \label{eq:ngek}
N_{\ge k} = O\left( \frac{n^4}{k^3} + \frac{n^6}{k^{11/2}}\log k + \frac{n^2}{k} \right) .
\end{equation}
(For the second term, put $Z:=m^3/N_{\ge k}$. If the second term in
the bound for $I(P,H)$ dominates, then we get $kN_{\ge k} =
O(m^{6/11}N_{\ge k}^{9/11}\log^{2/11} Z)$, or $k^{11/2} = O(Z\log
Z)$, from which one easily gets $Z = n^6/N_{\ge k} =
\Omega\left(k^{11/2}/\log k\right)$, and the bound follows.) Note
that for $3\le k<k_0$, we have the trivial bound $N_{\ge k} =
O(|P|^3) = O(n^6)$, as easily follows from the fact that none of our
planes intersect $\Pi$ in a line. Hence, (\ref{eq:ngek}) extends to
any $k\ge 3$, by suitably adjusting the constant of proportionality.

The analysis in Rudnev~\cite{Rud} uses the bound in (\ref{eq:ngek}) to derive an upper bound on
the number $Q$ of pairs of \emph{congruent pentuples} in $A$; that is, pairs of pentuples of elements of $A$
for which there exists a M\"obius transformation that maps the points in one pentuple to the respective points
in the second pentuple. Rudnev shows that $Q$ is upper bounded by
$$
Q = O\left( \sum_{j=2}^{\log n} N_{\ge 2^j}\cdot (2^j)^5 \right) .
$$
Using (\ref{eq:ngek}), we get
\begin{equation} \label{sump}
Q = O\left( \sum_{j=2}^{\log n} \left( 2^{2j}n^4 + \frac{jn^6}{2^{j/2}} + 2^{4j} n^2 \right) \right) = O(n^6) .
\end{equation}
In the original analysis, Rudnev has used, instead of (\ref{eq:ngek}), the weaker bound
$$
N_{\ge k} = O\left( \frac{n^6}{k^5} \right) .
$$
This bound follows from Solymosi and Tardos~\cite{SoTa}, and leads to the slightly inferior bound $Q= O(n^6\log n)$.

The rest of the analysis, from upper bounding $Q$ to obtaining a lower bound on the number of distinct cross-ratios,
is spelled out in \cite{Rud}. His original analysis yields the lower bound
$\Omega(n^{24/11}/\log^{6/11}n)$, and our improvement raises it to $\Omega(n^{24/11})$.

\medskip
\noindent{\bf Remark.}
For the purpose of showing that $P=O(n^6)$, we could equally use the more general bound in
Theorem~\ref{main}. The only difference is that the second term in (\ref{eq:ngek}) increases to
${\displaystyle O\left( \frac{n^6}{k^{11/2}}\log^{11\beta/2} k \right)}$, as is easily checked.
This does not affect the subsequent analysis, as the sum of the second terms in (\ref{sump})
is still the sum of a convergent series.

\paragraph{Acknowledgement.}
The authors would like to thank Misha Rudnev for valuable interaction and comments, especially
concerning the application in Section~\ref{sec:cross}.


\begin{thebibliography}{}

\bibitem{ANP+}
P. Agarwal, E. Nevo, J. Pach, R. Pinchasi, M. Sharir and S. Smorodinsky,
Lenses in arrangements of pseudocircles and their applications,
{\it J. ACM} 51 (2004), 139--186.

\bibitem{ApS}
R. Apfelbaum and M. Sharir,
Large bipartite graphs in incidence graphs of points and hyperplanes,
{\it SIAM J. Discrete Math.} 21 (2007), 707--725.

\bibitem{ArS}
B. Aronov and M. Sharir,
Cutting circles into pseudo-segments and improved bounds for incidences,
{\it Discrete Comput. Geom.} 28 (2002), 475--490.

\bibitem{BK}
P. Brass and Ch. Knauer,
On counting point-hyperplane incidences,
{\it Comput. Geom. Theory Appls.} 25 (2003), 13--20.

\bibitem{Chaz}
B. Chazelle,
Cuttings,
in {\it Handbook of Data Structures and Applications} (D.~P.~Mehta and S. Sahni, Eds.),
Chapman and Hall/CRC, 2004.


\bibitem{EGS}
H. Edelsbrunner, L. Guibas and M. Sharir,
The complexity of many cells in arrangements of planes and related problems,
{\it Discrete Comput. Geom.} 5 (1990), 197--216.

\bibitem{ET05}
Gy. Elekes and Cs.~D. T{\'{o}}th,
Incidences of not-too-degenerate hyperplanes,
in {\em Proc. 21st Annu. ACM Sympos. Comput. Geom.}, 2005, pages 16--21.

\bibitem{FPSSZ}
J. Fox, J. Pach, A. Sheffer, A. Suk, and J. Zahl, A semi-algebraic
version of Zarankiewicz's problem, {\it Journal of the European
Mathematical Society}, to appear.

\bibitem{Fu84}
W.~Fulton,
{\it Introduction to Intersection Theory in Algebraic Geometry},
Expository Lectures from the CBMS Regional Conference
Held at George Mason University, June 27--July 1, 1983, Vol. 54. AMS
Bookstore, 1984.


\bibitem{GK2}
L.\ Guth and N.\ H.\ Katz,
On the Erd{\H o}s distinct distances problem in the plane,
{\it Annals Math.} 181 (2015), 155--190. Also in arXiv:1011.4105.

\bibitem{Har}
C. G. A. Harnack,
\"Uber die Vieltheiligkeit der ebenen algebraischen Curven,
{\it Math. Ann.} 10 (1876), 189--199.


\bibitem{MT}
A. Marcus and G. Tardos,
Intersection reverse sequences and geometric applications,
{\it J. Combinat. Theory} Ser.~A 113 (2006), 675--691.

\bibitem{Rud}
M. Rudnev,
On distinct cross-ratios and related growth problems,
manuspcript, 2017.

\bibitem{SS3d}
M. Sharir and N. Solomon,
Incidences between points and lines in three dimensions,
in {\it Intuitive Geometry} (J. Pach, Ed.), to appear.
Also in {\it Proc. 31st Annu. Sympos. Computational Geometry} (2015), 553--568, and
in arXiv:1501.02544.


\bibitem{SZ}
M. Sharir and J. Zahl,
Cutting algebraic curves into pseudo-segments and applications,
{J. Combinat. Theory} Ser. A 150 (2017), 1--35.
Available online March 6, 2017.
Also in arXiv:1604.07877.


\bibitem{SoTa}
J. Solymosi and G. Tardos,
On the number of $k$-rich transformations,
{\it Proc. 23th ACM Sympos. Comput. Geom.}, 2007, 227--231.

\bibitem{Sze}
L. Sz\'ekely,
Crossing numbers and hard Erd\H os problems in discrete geometry,
{\it Combinat. Probab. Comput.} 6 (1997), 353--358.


\bibitem{Za}
J. Zahl,
An improved bound for the number of point-surface incidences in three dimensions,
{\it Contrib. Discrete Math.} 8(1) (2013), 100--121.

\end{thebibliography}
\end{document}